\documentclass[11pt]{amsart}

\usepackage{%
	amssymb%
	,euscript%
	,xspace%
	}

\newcommand{\ie}{i.e.,\xspace}

\newtheorem{thm}{Theorem}
\newtheorem{prop}{Proposition}
\newtheorem{lemma}{Lemma}
\newtheorem{cor}{Corollary}
\newtheorem{mainthm}{Main Theorem}

\theoremstyle{definition}
\newtheorem{defn}{Definition}
\newtheorem{prob}{Reduced Problem}

\theoremstyle{remark}
\newtheorem{ex}{Example}

\newcommand{\cal}{\EuScript}
\newcommand{\Hs}{\cal H}

\renewcommand{\leq}{\leqslant}
\renewcommand{\geq}{\geqslant}

\newcommand{\R}{\mathbb R}
\newcommand{\A}{\mathbb A}
\renewcommand{\H}{\mathbb H}
\newcommand{\X}{\mathbb X}
\newcommand{\C}{\mathbb C}
\newcommand{\T}{\mathbb T}

\newcommand{\norm}[1]{\left|#1\right|}

\newcommand{\innerh}[2]{\langle#1,#2\rangle_{\H}}
\newcommand{\innerx}[2]{\langle#1,#2\rangle_{\X}}

\newcommand{\dist}{\operatorname{dist}}
\newcommand{\compose}{\operatorname{\scriptstyle\circ}}

\newcommand{\Ric}{\operatorname{Ric}}
\newcommand{\tr}{\operatorname{tr}}

\newcommand{\p}{\varphi}
\newcommand{\D}{\partial}

\newcommand{\subsetnotequal}%
	{\subset\subset}

\newcommand{\bchi}{\boldsymbol{\chi}}
\newcommand{\bpsi}{\boldsymbol{\psi}}
\newcommand{\bF}{\boldsymbol{F}}
\newcommand{\bA}{\boldsymbol{A}}
\newcommand{\balpha}{\boldsymbol{\alpha}}
\newcommand{\bbeta}{\boldsymbol{\beta}}
\newcommand{\btheta}{\boldsymbol{\theta}}
\renewcommand{\v}{\boldsymbol{v}}

\newcommand{\bpsit}{\widetilde{\boldsymbol{\psi}}}
\newcommand{\bchit}{\widetilde{\boldsymbol{\chi}}}
\newcommand{\bvt}{\widetilde{\boldsymbol{v}}}

\newcommand{\pt}{\widetilde\varphi}
\newcommand{\vt}{\widetilde v}
\newcommand{\Sigmat}{\widetilde\Sigma}
\newcommand{\ut}{\widetilde u}

\newcommand{\ph}{\widehat\varphi}
\newcommand{\uh}{\widehat u}
\newcommand{\vh}{\widehat{\boldsymbol{v}}}

\newcommand{\Bb}{\overline B}

\newcommand{\aut}[1]{{\sc #1}}
\newcommand{\tit}[1]{{\em #1\/}}
\newcommand{\vol}[1]{{\bf #1}}
\newcommand{\yr}[1]{(#1)}
\newcommand{\pp}[2]{#1--#2}

\newcommand{\sm}{\setminus}

\newcommand{\step}[1]{\vspace{1em}\noindent{\bf Step #1.}}

\title[Harmonic Maps]{Harmonic
	Maps with Prescribed Singularities on Unbounded Domains}
\author{Gilbert Weinstein}
\address{Department of Mathematics, University of Alabama at Birmingham,
	Birmingham, Alabama 35205}
\thanks{This work was supported in part by NSF Grant DMS-9404523}
\email{weinstei@@math.uab.edu}
\date{\today}

\begin{document}

\maketitle

\begin{abstract}
The Einstein/Abelian-Yang-Mills Equations reduce in the stationary
and axially symmetric case to a harmonic map with prescribed singularities
$\p\colon\R^3\sm\Sigma\to\H^{k+1}_\C$ into the $(k+1)$-dimensional
complex hyperbolic space.
In this paper, we prove the existence and uniqueness of harmonic maps with
prescribed singularities $\p\colon\R^n\sm\Sigma\to\H$, where $\Sigma$ is an
unbounded smooth closed submanifold of $\R^n$ of codimension at least $2$,
and $\H$ is a real, complex, or quaternionic hyperbolic space.  As a
corollary, we prove the existence of solutions to the reduced stationary and
axially symmetric Einstein/Abelian-Yang-Mills Equations.
\end{abstract}

\section{Introduction}

It is well known that the Einstein/Maxwell Equations reduce in the stationary
and axially symmetric case to a harmonic map with prescribed singularities
$\p\colon\R^3\sm\Sigma\to\H^2_\C$ into the complex hyperbolic plane.
Here $\Sigma=\A\sm\bigcup_{j=1}^N I_j$ is the closed set obtained
by removing $N$ open intervals $I_j$ from the $z$-axis $\A\subset\R^3$.
More generally, the Einstein/Abelian-Yang-Mills
Equations with Abelian gauge group
$\T^k$ reduce in the same setting to a harmonic map with prescribed
singularities $\p\colon\R^3\sm\Sigma\to\H^{k+1}_\C$, where $\H^{k+1}_\C$ is
the $(k+1)$-dimensional complex hyperbolic space.
This model has $k$ `electromagnetic' fields which
interact weakly via the gravitational field.
The number $N$ corresponds to the number of black holes.

In~\cite{weinstein96}, we proved the existence and
uniqueness of harmonic maps with prescribed singularities
$\p\colon\R^n\sm\Sigma\to\H$, where $\H$ is any hyperbolic space (real,
complex, or quaternionic), and $\Sigma$ is any smooth closed submanifold of
$\R^n$ of codimension at least two, {\em provided\/} a harmonic map
$\pt\colon\R^n\sm\Sigmat\to\H$ was given,
prescribing the behavior of $\p$ at infinity in $\R^n\sm\Sigma$.
Here $\Sigmat$ denotes the union of the unbounded components of $\Sigma$.
The singular behavior of $\p$ near each component $\Sigma_j$ of
$\Sigma$ was prescribed
by a {\em $\Sigma_j$-singular map into $\gamma_j$\/}, that is a
harmonic map $\p_j\colon\R^n\sm\Sigma_j\to\H$
whose image lies within a geodesic $\gamma_j$ of $\H$, and
such that $\p_j(x)\to\gamma_j(+\infty)$ as
$x\to\Sigma_j$.  Such maps are easily constructed using harmonic functions.
The map $\p$ was required to be asymptotic to $\p_j$ near $\Sigma_j$,
and asymptotic to $\pt$ at infinity.
In particular, it followed that the map $\pt$ was to be asymptotic
near each unbounded component $\Sigma_j$ of $\Sigma$ to $\p_j$.
This map $\pt$ was to be considered
part of the data, and was to be found before the main theorem
in~\cite{weinstein96} could be applied.

For the harmonic map problem associated with the Einstein/Maxwell Equations,
the map $\pt$ could be obtained from an explicit solution, the Kerr-Newman
solution~\cite{newman65},
corresponding to $N=1$ black holes.  However, the existence of the harmonic
map $\pt$ in the more general situation of the main
theorem in~\cite{weinstein96} is not known.
In this paper, we replace the hypothesis on the existence of
the harmonic map $\pt$ with a considerably weaker hypothesis:
the existence of a map
$\pt$ which is only approximately harmonic at infinity, see the Main Theorem
in Section~\ref{sec:existence}.
Such maps are much easier to construct than harmonic maps.
As a corollary, we obtain the existence of solutions to the reduced
stationary and axially symmetric
Einstein/Abelian-Yang-Mills Equations corresponding to $N\geq1$ black holes,
\ie the existence of a harmonic map
with prescribed singularities $\p\colon\R^3\sm\Sigma\to\H^{k+1}_\C$,
where $\Sigma$ is as in the first paragraph of this introduction.  Of
course, when $k=1$, we recover our previous result from~\cite{weinstein96}.
Since we are interested in unbounded singular sets
$\Sigma\subset\R^n$ of codimension
at least $2$, we will assume throughout that $n\geq3$.

It is important to point out that, as before
in~\cite{weinstein94,weinstein96}, unless $N=1$, the
solutions of the Einstein/Abelian-Yang-Mills Equations constructed
in this way, are expected to have a conical singularity on some bounded
component of the axis.  This conical singularity can be interpreted as the
gravitational force between the black holes,
see~\cite{bach:weyl,weinstein94}.  However, when $N=1$, these spacetimes
should be regular, and asymptotically flat.  To establish this fact,
a regularity result as in~\cite{weinstein92,li:tian92} is
needed.  We will come back to this problem in a future paper.

In the next section, we describe briefly the setup
of the problem.  For more details, we refer the reader
to~\cite{weinstein95,weinstein96}.  In section~\ref{sec:existence},
we state our main existence result and give its proof.

\section{Set-Up and Preliminaries}	\label{sec:prem}

Let $n\geq3$, let $\Omega\subset\R^n$ be a smooth domain, and let
$\X$ be a Riemannian manifold of dimension $m$.
A {\em harmonic map\/} is a map $\p\colon\Omega\to\X$ which is for each
$\Omega'\subset\subset\Omega$
a critical point of the energy functional:
\[
	E_{\Omega'}(\p) = \int_{\Omega'} \norm{d\p}^2,
\]
where $\norm{d\p}^2=\sum_{j=1}^n \innerx{\D_j\p}{\D_j\p}$
is the energy density of $\p$,
and $\innerx{\cdot}{\cdot}$ represents the metric on $\X$.
If $\p$ is harmonic, it satisfies an elliptic system of nonlinear partial
differential equations
which can be written in local coordinates on $\X$ as:
\[
	\tau^a =
	\Delta \p^a  +
	\sum_{b,c=1}^m \sum_{j=1}^n \Gamma^a_{bc}(\p)
	\, \D_j\p^b\, \D_j\p^c =0,
\]
where $\Gamma^a_{bc}$ are the Christoffel symbols of $\X$.  Note that
for any map $\p$, the quantities
$\tau^a$ are the components of a cross section  $\tau(\p)$
of the pull-back bundle $\p^{-1} T\X$.
This vector field along $\p$ is called
the {\em tension\/} of $\p$, and its  magnitude, defined by
$\norm{\tau}^2 = \innerx\tau\tau$, is a real function on
$\Omega$.

Here we will be concerned only with harmonic maps into a classical globally
symmetric space of noncompact type and rank one, \ie a real, complex, or
quaternionic hyperbolic space.  We denote such a space by $\H$ and let $m$
be its real dimension.  The model we use for
$\H$ is $\R^m=\R\times\R^{m-1}$ where the first coordinate is a
{\em Busemann function\/}
on $\H$, see~\cite{eberlein,weinstein95}.  The metric of $\H$ then
takes the form:
\begin{equation}	\label{eq:hypmetric}
	ds^2 = du^2 + Q_p(dv),
\end{equation}
where for each $p\in\H$, $Q_p$ is a quadratic form on $\R^{m-1}$.

\begin{ex}
For example, when $\H$ is
the complex hyperbolic space $\H^{k+1}_\C$, we have $m=2k+2$, and we write
$p=(u,v,\bchi,\bpsi)$, where $\bchi,\bpsi\in\R^{k}$.
The quadratic form $Q_p$ is then given by:
\[
	Q_p(dv) = e^{4u} \left( dv + \bchi \cdot d\bpsi -
	\bpsi\cdot d\bchi\right)^2 + e^{2u} \left( d\bchi\cdot d\bchi +
	d\bpsi\cdot d\bpsi \right),
\]
where the dot represents the Euclidean inner product on $\R^k$,
see~\cite[Lemma 6]{weinstein95}.
We also give here for future reference the
components of the tension $\tau(\p)$ of a map $\p=(u,v,\bchi,\bpsi)$ into
$\H^{k+1}_\C$:
\begin{equation} \label{eq:tension}
\begin{aligned}
	\tau^u &= \Delta u - 2 e^{4u} \norm{\nabla v + \bchi\cdot\nabla\bpsi
	- \bpsi \cdot\nabla\bchi}^2 \\
	& \qquad\qquad\qquad\qquad\qquad\qquad {}
	- e^{2u} \left( \nabla\bchi\cdot\nabla\bchi +
	\nabla\bpsi\cdot\nabla\bpsi \right) \\
	\tau^v &= \Delta v + 4\, \nabla v \, \nabla u -
	2 e^{2u} \left( \bchi\cdot\nabla\bpsi + \bpsi\cdot\nabla\bchi
	\right) \, \nabla u \\
	&\qquad\qquad {} - 2^{2u} \left( \bchi\cdot\nabla\bchi +
	\bpsi\cdot\nabla\bpsi \right) \, \left( \nabla v
	+\bchi\cdot\nabla\psi - \bpsi\nabla\bchi \right) \\
	\tau^{\bchi} &=  \nabla\bchi + 2\, \nabla\bchi \, \nabla u - 2
	e^{2u} \nabla\bpsi \, \left( \nabla v + \bchi\cdot\nabla\bpsi -
	\bpsi\cdot\nabla\bchi \right)\\
	\tau^{\bpsi} &= \nabla\bpsi + 2\, \nabla\bpsi \, \nabla u + 2
	e^{2u} \nabla\bchi \, \left( \nabla v + \bchi\cdot\nabla\bpsi -
	\bpsi\cdot\nabla\bchi \right).
\end{aligned}
\end{equation}
\end{ex}

Harmonic maps have been studied extensively.  When the
target space is nonpositively curved, as is the case here,
finite energy harmonic maps are well
understood, see~\cite{eells:sampson,hamilton,schoen82,schoen83}.  In
particular, if $\p$ is any finite energy harmonic map then $\p$ is smooth.
In this paper, we deal with infinite energy harmonic maps into $\H$
with prescribed singularities along a submanifold $\Sigma$
of codimension at least two.  We now give a few definitions taken
from~\cite{weinstein95,weinstein96}.

\begin{defn}
Let $\Omega\subset\R^n$ be a smooth domain, and let
$\Sigma$ be a smooth closed submanifold of
$\R^n$ of codimension at least two.
Let $\gamma$ be a geodesic of $\H$.
We say that a harmonic map $\p\colon\Omega\sm\Sigma\to\H$ is a {\em
$\Sigma$-singular map into $\gamma$} if:
\begin{itemize}
	\item[(i)] $\p(\Omega\sm\Sigma)\subset\gamma(\R)$;
	\item[(ii)] $\p(x)\to\gamma(+\infty)$ as $x\to\Sigma$;
	\item[(iii)] There is a constant $\delta>0$
	such that $\norm{d\p(x)}^2 \geq \delta \dist(x,\Sigma)^{-2}$
	in some neighborhood of $\Sigma$.
\end{itemize}
\end{defn}

Note that if $u$ is a harmonic function on $\Omega\sm\Sigma$
then $\p=\gamma\compose u$ is a harmonic
map and $\norm{d\p}^2 = \norm{\nabla u}^2$.  Thus, the problem of
finding $\Sigma$-singular maps into $\gamma$ is
simply a problem in potential theory.

\begin{ex}	\label{ex:sigmasing}
Let $\Omega\subset\R^n$ and let $\Sigma\subset\Omega$ be a smooth closed
submanifold of codimension at least two.  Let $u_0$ be the potential
of a charge distribution which is positive and bounded away from
zero on $\Sigma$, and let $\v_0\in\R^{m-1}$ be constant.
Then, using the coordinate system for $\H$ described above in which the
metric takes the form~\eqref{eq:hypmetric}, we see that
$\p=(u_0,\v_0)$ is a $\Sigma$-singular map
into the geodesic $\gamma(t)=(t,\v_0)$ of $\H$.  Such maps will be used
to prescribe the singular behavior near $\Sigma$.
\end{ex}

The next definition generalizes the concept of asymptotic geodesics, see for
example~\cite{eberlein}.

\begin{defn}
Let $\p,\p'\colon\Omega\sm\Sigma\to\H$ be maps, and let
$\Sigma'\subset\Sigma$.  We say that $\p$ and $\p'$ are
{\em asymptotic near $\Sigma'$\/} if there is a neighborhood
$\Omega'$ of $\Sigma'$ such that $\dist(\p,\p')\in
L^\infty(\Omega'\sm\Sigma')$.
\end{defn}

To ensure uniqueness, we will use the following notion
for being asymptotic at infinity.

\begin{defn}
Let $\p,\p'\colon\R^n\sm\Sigma\to\H$ be maps.
We say that $\p$ and $\p'$ are {\em asymptotic at infinity\/},
if $\dist(\p,\p')\to0$ as $x\to\infty$ in $\R^n\sm\Sigma$.
\end{defn}

Consider now the Einstein/Abelian-Yang-Mills Equations
with an Abelian gauge group $\T^k$:
\begin{equation}	\label{eq:einstein}
\begin{gathered}
	\Ric_g - \frac{1}{2} \, R_g \, g = 2 T_{\bF}, \\
	\bF = d\bA, \\
	d\ast\bF =0, \\
	T_{\bF}(X,Y) = \tr_g (i_X \bF\cdot i_Y \bF
	+ i_X\ast \bF\cdot i_Y\ast \bF).
\end{gathered}
\end{equation}
Here $g$ is a Lorentzian metric on a 4-dimensional manifold $M$,
$\Ric_g$ is the Ricci curvature of $g$, $R_g$ its scalar curvature, $\bA$
is an $\R^k$-valued one-form on $M$,
$d$ is the exterior derivative, $i_X$ inner
multiplication by the vector $X$,
$\tr_g$ is the trace with respect to the metric $g$,
and the dot represents a
Euclidean inner product on $\R^k$.  We impose the condition that the
solution is asymptotically flat, stationary and axially symmetric, globally
hyperbolic, and that its domain of outer communications is simply connected.
Furthermore, we will assume that $(M,g)$ has an event horizon with
$N\geq 1$ connected components, each of
which is nondegenerate, see~\cite{weinstein96} for the definitions.

Following step by step the reduction
in~\cite[Section 2]{weinstein96}, one sees
that the equations reduce to a harmonic map with prescribed singularities
into $\H^{k+1}_\C$.  For this, one only needs to note
that the one-forms $\balpha=i_\xi \bF$, and
$\bbeta=i_\xi\ast \bF$ are now $\R^k$-valued.
Here $\xi$ is the generator of the axial symmetry normalized so that its
orbits have parameter length $2\pi$.
Thus the potentials $\bchi$, and $\bpsi$, given by
$d\bchi=\balpha$, and $d\bpsi=\bbeta$ are also $\R^k$-valued.
{}From here on, all the calculations go through,
with all products between $\R^k$-valued objects being taken in the sense
of the Euclidean inner product on $\R^k$.
In particular, the {\em twist\/} of $\xi$
is given by $\omega=\ast(d\xi\wedge\xi)$, where $\ast$ is the Hodge
star operator, and the {\em twist potential\/} $v$ is defined by
$2\, dv = \omega - 2 (\bchi\cdot d\bpsi - \bpsi\cdot d\bchi)$.
Define $u=-\log\norm{\xi}$ on $M'=\{g(\xi,\xi)>0\}$, and let $Q$ be the
quotient of $M'$ by its group $\R\times SO(2)$ of isometries.    Then
$Q\times SO(2)$ is diffeomorphic to $\R^3\sm\Sigma$, where
$\Sigma=\A\sm\bigcup_{j=1}^N I_j$, $I_j$ are open intervals of the
$z$-axis $\A\subset\R^3$, and $N$ is the number of connected components of
the event horizon.  Let $\tau$ be another generator of the group of
isometries, normalized so that $\norm{\tau}^2\to-1$ at spacelike infinity in
$(M,g)$, and let $\rho$ be the area element of the orbit, \ie $\rho^2 = -
\norm{\xi\wedge\tau}^2$.  Then $\Delta\rho=0$, hence $\rho$ is a harmonic
coordinate on $Q$, and we can pick a conjugate
harmonic coordinate $z$, so that the
metric of $Q$ is conformal to $d\rho^2 + dz^2$.  We now put on $Q\times
SO(2)=\R^3\sm\Sigma$ the flat metric $d\rho^2+dz^2 + \rho^2 d\phi^2$.
Let $\p=(u,v,\bchi,\bpsi)$, then $\p$ is
an axially symmetric harmonic map into $\H^{k+1}_\C$ defined on
$\R^3\sm\Sigma$.  Thus, we obtain the following result analogous to
Theorem~1 in~\cite{weinstein96}.

\begin{thm}	\label{thm:ernst}
Let $(M,g,\bA)$ be
a solution of the Einstein/Abelian-Yang-Mills Equations with Abelian
gauge group $\T^k$, which is stationary and axially symmetric.
Assume that $M$ is simply connected, that $(M,g)$ is the domain of outer
communications of an asymptotically flat, globally
hyperbolic spacetime, and that the event horizon in
$(M,g)$ has $N\geq1$ connected components each of which
is nondegenerate.  Define $\Sigma\subset\R^3$ and
$\p\colon\R^3\sm\Sigma\to\H^{k+1}_\C$ as above.  Then
$\p$ is an axially symmetric harmonic map.
\end{thm}

Furthermore the axis regularity conditions imply that the map $\p$ is
asymptotic near each $\Sigma_j$ to a $\Sigma_j$-singular harmonic map
$\p_j$ into a geodesic $\gamma_j$ of $\H^{k+1}_\C$.
To prescribe the behavior at infinity, we now construct a map
$\pt\colon\R^3\sm\Sigma\to\H^{k+1}_\C$ which is nearly harmonic at infinity.
We note that essentially the same construction will
give maps from $\R^n\sm\Sigma$ into quaternionic hyperbolic space, provided
$\Sigma$ is of codimension 2, and each component of
$\Sigma$ is contained within some cone centered at the origin.

\begin{ex}	\label{ex:quasi}
Let $\Sigma\subset\R^3$ be as above.  For
each $1\leq j\leq N+1$, let $\p_j\colon\R^3\sm\Sigma\to\H^{k+1}_\C$ be a
given $\Sigma_j$-singular map into a geodesic $\gamma_j$ of
$\H^{k+1}_\C$.  We assume, without loss of generality,
that all the geodesics $\gamma_j$ have a common
initial point $\gamma_j(-\infty)$ on $\D\H^{k+1}_\C$, so that the maps
$\p_j$ can all be written as
in Example~\ref{ex:sigmasing}, \ie $\p_j=(u_0,v_j,\bchi_j,\bpsi_j)$ where
$v_j\in\R$, and $\bchi_j,\bpsi_j\in\R^k$ are constants.  We also
assume that $u_0$ is the potential of a
uniform unit charge on $\Sigma$.  We can then normalize $u_0$ so that
$u_0\to-\log\rho$ as $r\to\infty$ in $\R^3\sm\Sigma$.  We seek a map
$\pt\colon\R^3\sm\Sigma\to\H^{k+1}_\C$
which coincides in a neighborhood of each $\Sigma_j$
with $\p_j$, and which is nearly harmonic at infinity.
More precisely, we will
obtain a map $\pt$ for which there is a constant $c>0$ such that
\begin{equation}	\label{eq:quasi}
	\norm{\tau(\pt)}\leq c (1+r^2)^{-2},
\end{equation}
where $r$ is the distance to the origin in $\R^3$.
To achieve this, we simply set $\ut=u_0$, and
we pick functions $\vt,\bchit,\bpsit$ which depend only
on the polar angle outside a sufficiently large ball, and
which take on the appropriate constant values $v_j,\bchi_j,\bpsi_j$,
near each $\Sigma_j$.  It is now a simple calculation
using~\eqref{eq:tension} to check that
the map $\pt=(\ut,\vt,\bchit,\bpsit)$ satisfies~\eqref{eq:quasi}.
\end{ex}

We can now pose the following Reduced Problem for the stationary and axially
symmetric Einstein/Abelian-Yang-Mills Equations.

\begin{prob}
Let $\Sigma=\bigcup_{j=1}^{N+1} \Sigma_j =
\A\sm\bigcup_{j=1}^N I_j$, where $\A$ is the $z$-axis in $\R^3$,
$I_j$ are $N$ open intervals, and $\Sigma_j$ are the connected
components of $\Sigma$.  Let
$u_0$ be the potential of a uniform unit charge
distribution on $\Sigma$.  For each $1\leq j\leq N+1$, let
$v_j\in\R$ and $\bpsi_j\in\R^k$ be constants, and let
$\p_j\colon\R^3\sm\Sigma\to\H^{k+1}_\C$ be the $\Sigma$-singular map given,
as in Example~\ref{ex:sigmasing},
by $\p_j=(u_0,v_j,\boldsymbol{0},\bpsi_j)$.
Construct the map $\pt\colon\R^3\sm\Sigma\to\H^{k+1}_\C$
as in Example~\ref{ex:quasi}.
The Reduced Problem is then: to find an axially symmetric harmonic map
$\p\colon\R^3\sm\Sigma\to\H^{k+1}_\C$ which is asymptotic to $\pt$ near
each $\Sigma_j$, and at infinity.
\end{prob}

Let $\Sigma'=\Sigma_1\cup\Sigma_{N+1}$, and
denote by $\pt'\colon\R^3\sm\Sigma'\to\H^2_\C$
the harmonic map obtained from the Kerr-Newman
solution~\cite{carter73},
which was used in~\cite{weinstein96} to prescribe the behavior of
$\p$ at infinity.
Note that $\pt$ and $\pt'$ are asymptotic at infinity.  Thus,
using $\pt'$ is equivalent in this case to using $\pt$.

In the next section, we prove a theorem which implies that the Reduced
Problem has a unique solution for each value of the $(N+1)(k+1)$
parameters.  We note that one of the constants $v_j$ and one of the
constants $\bpsi_j$ can be set to zero using an isometry of $\H^{k+1}_\C$,
so that we have left $N(k+1)$ parameters.
These correspond to $Nk$ charges and $N$
angular momenta.  The $N$ distances and the $N$ masses are fixed by
the choice of the intervals $I_j$.

There is a partial converse to Theorem~\ref{thm:ernst}.
Let $\p=(u,v,\bchi,\bpsi)
\colon\R^3\sm\Sigma\to\H^{k+1}_\C$ be a solution of the Reduced
Problem.  Let $(\rho,\phi,z)$ be cylindrical coordinates on $\R^3$.
If $\omega = 2(dv + \bchi\cdot d\bpsi + \bpsi\cdot d\bchi)$, then
it follows from the harmonic map equations that the two-form
$e^{4u}i_\xi\ast\omega$ is closed,
and hence there is a function $w$ such that
$dw=e^{4u}i_\xi\ast\omega$.  Here $\ast$ is the Hodge star operator of the
Euclidean metric on $\R^3$.  Similarly, one checks that there is a function
$\lambda$ such that
\begin{multline*}
	d\lambda = du +
	\rho\left[ u_\rho^2 - u_z^2 + \frac{1}{4} e^{4u}(\omega_\rho^2 -
	\omega_z^2) \right. \\
	\left. \phantom{\frac14} + e^{2u} (\bchi_\rho\cdot\bchi_\rho -
	\bchi_z\cdot\bchi_z + \bpsi_\rho\cdot\bpsi_\rho -
	\bpsi_z\cdot\bpsi_z) \right] d\rho \\[.5ex]
	 + 2 \rho \left[
	u_\rho u_z + \frac{1}{4} e^{4u} \omega_\rho \, \omega_z +
	e^{2u} (\bchi_\rho \cdot \bchi_z
	+ \bpsi_\rho \cdot \bpsi_z) \right] dz.
\end{multline*}
Also, there is an $\R^k$-valued
function $\btheta$ such that $d\btheta=e^{2u}i_\xi\ast
d\bpsi - w\, d\bchi$.  As before, the forms $d\lambda$ and $d\btheta$ are
closed thanks to the harmonic map equations.
We now define the metric $g$ on
$M'=\R\times(\R^3\sm\A)$ by:
\[
	ds^2 = -\rho^2 e^{2u}\, dt^2 + e^{-2u} (d\phi - w\, dt)^2 +
	e^{2\lambda} ( d\rho^2 + dz^2),
\]
and the $\R^k$-valued one-form $\bA$ by:
\[
	\bA = -( \bchi\, d\phi + \btheta\, dt).
\]
We have:

\begin{thm}	\label{thm}
Let $\p=(u,v,\bchi,\bpsi)
\colon\R^3\sm\Sigma\to\H^{k+1}_\C$ be a solution of the Reduced
Problem, and define $M$, $g$, and $\bA$ as above.  Then $(M,g,\bA)$ is a
solution of the Einstein/Abelian-Yang-Mills Equations.
\end{thm}

The intervals $I_j\subset\A$ correspond to an event horizon in $(M,g)$.
However, as we mentioned in the introduction, unless $N=1$,
it is expected that the metric $g$ cannot be extended across
all of $\Sigma$.  The
obstruction is a conical singularity on some bounded components of $\Sigma$.
This singularity can be related to the force between the black holes,
see~\cite{bach:weyl,weinstein94}.
Nevertheless, $(M,g)$ should be asymptotically flat, since a conical
singularity cannot occur on the two unbounded components of $\Sigma$.  To
establish asymptotic flatness, the first step is a `regularity' result for
$\p$ across $\Sigma$ as in~\cite{weinstein92,li:tian92}.  This will be
addressed in a future paper.

\section{Results and Proofs}	\label{sec:existence}

In this section, we prove the following theorem which is a generalization of
the Main Theorem in~\cite{weinstein96}.

\begin{mainthm}
Let $\Sigma$ be a closed smooth submanifold of $\R^n$ of codimension
at least two,
and let $\Sigma_j$ denote its connected components.
For each $1\leq j\leq N$, let $\p_j\colon\R^n\sm\Sigma\to\H$ be a
$\Sigma$-singular harmonic map into some geodesic $\gamma_j$ of the
space $\H$.  Suppose there is a map $\pt\colon\R^n\sm\Sigma\to\H$
which coincides in a neighborhood of each $\Sigma_j$
with $\p_j$, and for which there is a constant $c>0$ such that
\begin{equation}	\label{eq:quasi3/2}
	\norm{\tau(\pt)}\leq c (1+r^2)^{-3/2}.
\end{equation}
Then, there is a unique harmonic map $\p\colon\R^n\sm\Sigma\to\H$ which is
asymptotic to $\pt$ near each $\Sigma_j$, and at infinity.
\end{mainthm}

In view of Example~\ref{ex:quasi},
we obtain as an immediate corollary the existence and uniqueness of
solutions to the Reduced Problem.

\begin{cor}
The Reduced Problem for the stationary and axially symmetric
Einstein/Abelian-Yang-Mills Equations has a unique solution.
\end{cor}

\noindent{\em Proof of the Main Theorem.}
The existence proof is divided, as in~\cite{weinstein96}, into
three steps:
\begin{description}
\item[Step 1]
for each sufficiently large ball $B_R\subset\R^n$, there exists a
map $\ph_R\colon B_R\sm\Sigma\to\H$ such that
$\ph_R$ is asymptotic to $\pt$
near each $\Sigma_j$, and such that $\ph_R=\pt$ on $\D B_R$.
\item[Step 2]
there is a constant $C$ independent of $R$ such that
$\dist(\ph_R,\pt)\leq C$ in $B_R\sm\Sigma$.
\item[Step 3]
there is a sequence $R_i\to\infty$ for which $\ph_{R_i}$
converge uniformly on compact subsets of $\R^n\sm\Sigma$
to a harmonic map $\ph\colon\R^n\sm\Sigma\to\H$
which is asymptotic to $\pt$
near each $\Sigma_j$, and at infinity.
\end{description}

\step{1}
This follows from Proposition~1 in~\cite{weinstein96}.  Although,
we required there that $\pt$ was harmonic, this was never used in
the proof of Proposition~1.
We also required that $\dist(\pt(x),\gamma_j)\to0$ as
$x\to\Sigma_j$, but that is clearly satisfied
here since $\pt$ agrees with $\p_j$
near $\Sigma_j$.  In fact, the existence of $\ph_R$ is proved
using a variational approach, and thus one also obtains in the
course of the proof an energy estimate which is used in Step 3.

We may assume, without loss of generality, that all the geodesics $\gamma_j$
have a common initial point $p=\gamma_j(-\infty)\in\D\H$.
We use on $\H$ the coordinate system $(u,\v)\colon\H\to\R\times\R^{m-1}$
mentioned in Section~\ref{sec:prem}, where $u$ is the Busemann function
associated with the point $p$:
\[
	u(q) = \lim_{t\to-\infty}\bigl(\dist\bigl(q,\gamma_1(t)\bigr) +
	t\bigr),
\]
see~\cite[Lemma 6]{weinstein95}.  The metric then takes the
form~\eqref{eq:hypmetric}.  Write $\pt=(\ut,\bvt)$, and $\p_j=(u_0,\v_j)$,
where $u_0$ is a harmonic function with $u_0(x)\to\infty$ as $x\to\Sigma$.
Here, we have assumed that all the maps $\p_j$ have the same harmonic
function $u_0$ as their $u$ coordinates, but that can easily be arranged.

Let $H_1(B_R)$ be the Sobolev space of functions $u$ on $B_R$
such that $u,\nabla u\in L^2(B_R)$, and denote the closure
of $C^\infty_0(B_R)$ in that space by $H_{1,0}(B_R)$.
Let $H_1^{\pt}(B_R,\R^{m-1})$ be the weighted Sobolev space of
$\R^{m-1}$-valued functions $\v$ on $B_R$ such that
\[
	\int_{B_R} \left\{\norm{\v}^2 + Q_{\pt}(\nabla\v)\right\} < \infty,
\]
and let the closure of $C^\infty_0(B_R\sm\Sigma,\R^{m-1})$
in that space be $H_{1,0}^{\pt}(B_R,\R^{m-1})$.
Now, define the space $\Hs(B_R)$ of maps $\p=(u,\v)$ on $B_R$ satisfying
\[
	\begin{cases}
	u-\ut \in H_1(\Omega) \\
	\v-\bvt \in H_{1,0}^{\pt}(B_R;\R^{m-1}) \\
	\dist(\p,\pt) \in L^\infty(B_R),
	\end{cases}
\]
and let $\Hs_C(B_R)$ be the space of maps $\p\in\Hs(B_R)$ for which
$\dist(\p,\pt)\leq C$.
Define the renormalized energy of $\p=(u,\v)\in\Hs_C(B_R)$ by:
\[
	F_R(\p) = \int_{B_R} \left\{
	\norm{\nabla(u-u_0)}^2 + Q_{\p}(\nabla\v)
	\right\}.
\]
The solution is found by minimizing $F_R(\p)$ over $\p\in\Hs(B_R)$.  In
fact, an almost standard
direct variational argument shows that a minimizer of $F_R$ exists
in $\Hs_C(B_R)$ for any $C>0$.
The main difficulty is in proving that for some
$C=C_R$ large enough,
depending on $R$,
the infimum of $F_R$ over $\Hs_C(B_R)$ is the same as over $\Hs(B_R)$.

\begin{prop}	\label{prop:BR}
Let $R>0$ be large enough so that each bounded component of $\Sigma$ is
contained in the ball $B_R$.  Then there is a unique harmonic map
$\ph_R\colon
B_R\sm\Sigma\to\H$ such that $\ph_R$ is asymptotic to $\pt$ near each
$\Sigma_j$, and such that $\ph_R=\pt$ on $\D B_R$.  Furthermore,
$\ph_R\in\Hs(B_R)$.
\end{prop}

\step{2}
The main difference with~\cite{weinstein96} lies in this step.
We must now show that the solutions $\ph_R$
obtained in Step 1 belong to the space
$\Hs_C(B_R)$ with the constant $C$ independent of $R$.  For this purpose,
we will need two simple lemmas.

\begin{lemma}	\label{lemma:dist}
Let $\p_1,\p_2\colon\Omega\to\H$ be smooth maps,
and let $\rho=\dist(\p_1,\p_2)$.  Then, we have
\begin{equation}
	\label{eq:rho}
	\Delta\sqrt{1+\rho^2}
	\geq - \biggl( \norm{\tau(\p_1)} +
	\norm{\tau(\p_2)} \biggr)
\end{equation}
\end{lemma}

\noindent{Proof.}
Let $G_p(q)$ denote the gradient of the function
$q\mapsto\dist(p,q)$ on $\H$ evaluated at $q\ne p$,
and let $H$ denote the Hessian of the function $(p,q)\mapsto\dist(p,q)$ on
$\H\times\H$ as a quadratic form on $T(\H\times \H)$ away from the diagonal
$\{(p,p)\in\H\times\H\}$.  The proof of the lemma relies
on the following identity:
\[
	\Delta\rho = \innerh{G_{\p_1}(\p_2)}{\tau(\p_2)}
			+ \innerh{G_{\p_2}(\p_1)}{\tau(\p_1)}
			+ H(d\p_1+d\p_2,d\p_1+d\p_2),
\]
which holds wherever $\rho\ne0$.  Equation~\eqref{eq:rho} now follows from
the fact that  $\norm{G_p(q)}=1$ for all $q\ne p$, and the fact
that $H\geq0$, due to the negative curvature of $\H$,
see~\cite[p.~368]{schoen:yau79}.
\qed

The next lemma is a simple generalization of the maximum principle for
subharmonic functions.

\begin{lemma}	\label{lemma:maxprinc}
Suppose $\sigma\in C^\infty(\Bb_R\sm\Sigma)$ is bounded and satisfies:
\begin{alignat*}{2}
	\Delta \sigma &\geq - c(1+r^2)^{-3/2}, &
		\quad\text{on $B_R\sm\Sigma$} & \\
	\sigma &= 0, &\quad\text{on $\D B_R$}. &
\end{alignat*}
Then, there holds:
\begin{equation}	\label{eq:bound}
	\sigma \leq \frac{c}{n-2}, \qquad\text{on $B_R$}.
\end{equation}
\end{lemma}

\noindent{pf}
Clearly, it suffices to prove the lemma with $c=1$.  Let $\nu$ be the
solution of
\begin{align*}
	\Delta \nu &= (1+r^2)^{-3/2} \\
	\nu &\to 0, \quad\text{as $x\to\infty$}.
\end{align*}
Note that such a solution exists, is radially symmetric, negative, bounded,
and satisfies
\begin{equation}	\label{eq:nu}
	\inf_{B_R} \nu = \nu(0) = -\frac{1}{n-2}.
\end{equation}
Now, we have $\Delta(\sigma+\nu)\geq0$ on $B_R\sm\Sigma$, and $\sigma+\nu$
is bounded.  Thus, we can show, using cut-off functions near
$\Sigma$, as in~\cite[Lemma 7]{weinstein96}, that $\sigma+\nu$ is weakly
subharmonic on all of $B_R$, with finite Dirichlet integral over
any ball $B_{R'}\subsetnotequal B_R$.
Thus, using the maximum principle, we obtain:
\begin{equation}	\label{eq:maxprinc}
	\sigma + \nu \leq \nu(R) \leq 0.
\end{equation}
In view of~\eqref{eq:nu}, this implies~\eqref{eq:bound}.
\qed

Step 2 now follows easily, for if $\rho=\dist(\ph_R,\pt)$,
and $\sigma=\sqrt{1+\rho^2}-1$ then $\sigma$ is bounded, and
$\sigma=0$ on $\D B_R$.  Using~\eqref{eq:quasi3/2}, it follows from
Lemma~\ref{lemma:dist} that $\sigma$ satisfies:
\[
	\Delta \sigma \geq -c(1+r^2)^{3/2}.
\]
Thus, by Lemma~\ref{lemma:maxprinc}, there is a constant $C$,
independent of $R$, such that $\rho\leq C$.
In Step 3, we will also use the more precise estimate~\eqref{eq:maxprinc}.

\step{3}
Again, the proof of this step is very similar to
the one in~\cite{weinstein96}, and we only give a sketch.
By Steps 1 and 2, we have for each $R$ large enough,
a harmonic map $\ph_R=(\uh_R,\vh_R)\in\Hs_C(B_R)$.  Now fix $R_0$, and for
$R>R_0$ consider the maps $\ph_{R}$ restricted to $B_{R_0}$.  It
is not difficult to check, as in~\cite{weinstein96}, that
$F_{R_0}(\ph_R)$ is uniformly bounded, \ie
bounded independently of $R$.  In fact, we
have $\Delta(\uh_R-u_0)\geq0$, and
\[
	\norm{\uh_R-u_0}\leq \dist(\ph_R,\pt)+\norm{\ut-u_0} \leq C +
	\sup_{B_{R_0}} \norm{\ut-u_0},
\]
on $B_{R_0}\sm\Sigma$.
Thus, $\uh_R-u_0$ are uniformly bounded subharmonic functions on
$B_{R_0}\sm\Sigma$, and the argument used in the proof
of~\cite[Lemma 7]{weinstein96}, shows that the Dirichlet integrals of
$\uh_R-u_0$ over $B_{R_0}$
are uniformly bounded.  The bound on the second term in $F_{R_0}(\ph_R)$
now follows from the inequality:
\[
	2 Q_{\ph_R}(\nabla\vh_R) \leq \Delta(\uh_R-u_0).
\]
Therefore, we can find a sequence $R_i\to\infty$ such that
$\ph_{R_i}$ converges pointwise a.e.\ in $B_R$.
By a standard diagonal argument, we can assume
that the same sequence works for
all $R_0$.  It is not difficult to see that the pointwise limit $\p$ is a
harmonic map on $\R^n\sm\Sigma$.
Indeed if $\Omega\subset\subset\R^n\sm\Sigma$, then for $i$
large enough, $\ph_{R_i}|_{\Omega}$ is a family of smooth harmonic maps with
uniformly bounded energy which maps into a fixed compact set of $\H$.
Standard harmonic map theory now implies uniform bounds in
$C^{2,\alpha}(\Omega)$, hence some subsequence converges uniformly, by
necessity to $\p$, together with two derivatives.  Thus, $\p$ is harmonic.
Clearly, $\dist(\p,\pt)\leq C$ hence $\p$ is asymptotic to $\pt$ near each
$\Sigma_j$.  It remains to see that $\p$ is asymptotic to $\pt$ at
infinity.  This follows from~\eqref{eq:maxprinc} applied to
\[
	\sigma_R=\sqrt{1+\dist(\ph_R,\pt)^2}-1.
\]
In the limit we obtain
\[
	\sigma=\sqrt{1+\dist(\p,\pt)^2}-1\leq -c\, \nu,
\]
which implies that
$\dist(\p,\pt)\to0$ as $x\to\infty$ in $\R^n\sm\Sigma$.

The proof of the uniqueness statement is unchanged from~\cite{weinstein96}.
\qed

\end{document}